\documentclass{article}
\usepackage{amssymb,amsmath,amsthm,graphicx}
\usepackage[all,color]{xy}

\def\qed{\hfill {\hbox{${\vcenter{\vbox{               
   \hrule height 0.4pt\hbox{\vrule width 0.4pt height 6pt
   \kern5pt\vrule width 0.4pt}\hrule height 0.4pt}}}$}}}

\newtheorem{theorem}{Theorem}

\newtheorem{proposition}[theorem]{Proposition}

\theoremstyle{definition}
\newtheorem{example}{Example}
\newtheorem{definition}{Definition}
\newtheorem{remark}{Remark}

\date{}

\title{\Large \textbf{Psybrackets, Pseudoknots and Singular Knots}}

\author{
Suhyeon Jeong\footnote{Email: j00399501303@pusan.ac.kr. Supported by Basic Science Research Program through the National Research Foundation of Korea (NRF), funded by the Ministry of Education, Science and Technology (NRF-2019R1F1A1060205).} \and
Jieon Kim\footnote{Email: jieonkim7@gmail.com. Supported by Young Researchers Program through the National Research Foundation of Korea (NRF), funded by the Ministry of Education, Science and Technology (NRF-2018R1C1B6007021).} \and
Sam Nelson\footnote{Email: Sam.Nelson@cmc.edu. Partially supported by Simons Foundation collaboration grant 316709.}}

\begin{document}
\maketitle

\begin{abstract}
We introduce algebraic structures known as \textit{psybrackets} and use them 
to define invariants of pseudoknots and singular knots and links. Psybrackets
are Niebrzydowski tribrackets with additional structure inspired by the
Reidemeister moves for pseudoknots and singular knots. Examples and computations
are provided.
\end{abstract}

\parbox{5.5in} {\textsc{Keywords:} Pseudoknots, Singular knots, Psybrackets, 
Niebrzydowski tribrackets, ternary quasigroups

\smallskip

}

\section{\large\textbf{Introduction}}\label{I}

In \cite{MN,MN2}
algebraic structures known as \textit{knot-theoretic ternary quasigroups}
were introduced and investigated. With a notational change, these have 
been studied by the third listed author and collaborators as 
\textit{Niebrzydowski tribrackets} in papers such as \cite{GNT,NOO,NP,NP2} 
and used to define invariants of classical knots and links, virtual links 
and handlebody-links. Related objects known as \textit{biquasiles} have 
been investigated by the first and third authors in \cite{KN} and by the
third author and collaborators in \cite{CNS,DSN2}.

\textit{Pseudoknots} arose in biology as a way of dealing with knotted 
objects with only partial information about the crossings; see e.g.
\cite{CDDK,EP,LR,QR} etc. The mathematical formulation in \cite{H,HHJJMR,HJ}
defines pseudoknots and pseudolinks combinatorially as equivalence classes of
\textit{pseudoknot diagrams}, i.e., knot diagrams with ordinary classical 
crossings together with \textit{precrossings} in which it is unknown which
strand is on top, under the equivalence relation determined by the
\textit{pseudoknot Reidemeister moves}. 

\textit{Singular knots} are rigid vertex isotopy classes of 4-valent spatial 
graphs. We can think of singular knots and links as knots and links in which 
some strands are fused together at vertices known as \textit{singular 
crossings}. In particular, the cyclic ordering of the edges around each 
singular crossing is fixed. 

Identifying singular crossings with precrossings, the singular 
Reidemeister moves form a subset of the pseudoknot Reidemeister moves; 
combinatorially, the two classes of objects differ only by a single move.

In \cite{NOR}, together with two collaborators the third listed author 
exploited the similarity of the Reidemeister moves for pseudoknots
and singular knots to
introduce \textit{psyquandles}, algebraic coloring structures for 
pseudoknots and singular knots extending the notion of \textit{biquandle 
colorings} from the world of classical knots and links. Finite biquandles
give rise to integer-valued \textit{counting invariants}, which can be 
enhanced in various ways to define new stronger invariants.

In this paper we apply the idea of Niebrzydowski tribrackets to the case of
pseudoknots and singular knots, defining \textit{psybrackets} analogously
to the way psyquandles extend biquandles. The paper is organized as follows.
In Section \ref{P} we review the basics of pseudoknots and singular knots.
In Section \ref{PSY} we define psybrackets and provide some examples. In 
Section \ref{PCI} we define the psybracket counting invariant and provide some
computational examples to explore the power of the new invariants. We 
conclude in Section \ref{Q} with some questions for future research.

\section{\large\textbf{Pseudoknots and Singular Knots}}\label{P}

In this section we review the basics of \textit{pseudoknots} and 
\textit{singular knots}; the remainder of the paper will concern
algebraic structures from which we will derive invariants of 
pseudoknots and singular knots.

\begin{definition}
An \textit{oriented pseudoknot diagram} has \textit{positive and negative
classical crossings} but also \textit{precrossings}.
\[\includegraphics{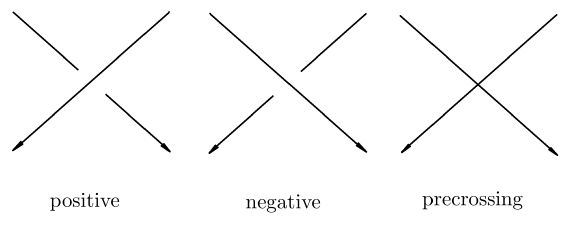}\]
Replacing a precrossing with a classical crossing is known as 
\textit{resolution}.
\end{definition}

Precrossings represent classical crossings for which it is unknown 
which strand passes over and which passes under. We may regard a precrossing
as a linear combination of both crossings with a scalar coefficient 
of $\frac{1}{2}$ for each; extending linearly, we may regard a pseudoknot
diagram as a linear combination of its resolutions. Interpreting the
scalar weights as probabilities, we obtain from a pseudoknot its
\textit{wereset} or \textit{weighted resolution set}, a  
discrete probability distribution whose events are the classical knots
obtained by resolving all precrossings, with probabilities 
given by the scalar coefficients.

\begin{example}
The pseudolink 
\raisebox{-0.4in}{\includegraphics{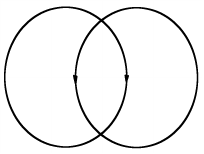}} has wereset
\[\left\{
\frac{1}{4}\raisebox{-0.4in}{\includegraphics{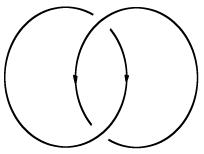}},
\frac{1}{4}\raisebox{-0.4in}{\includegraphics{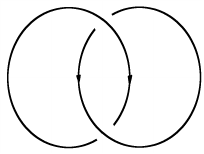}},
\frac{1}{2}\raisebox{-0.4in}{\includegraphics{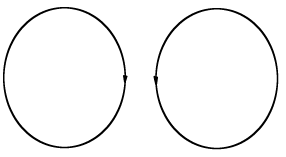}}
\right\}.\]
\end{example}

Singular knots and links are 4-valent spatial graphs considered up to
\textit{rigid vertex isotopy}, where the cyclic ordering of the edges 
entering a vertex is fixed. Such a rigid vertex is called a \textit{singular 
crossing}; we can imagine singular crossings as points where the knot becomes 
stuck to itself (transversely, not tangentially).

Identifying precrossings with singular crossings, the Reidemeister moves for 
pseudolinks and singular links are the same except for one move: precrossings
can be introduced or removed via a Reidemeister I type move, while singular
crossings cannot. More precisely, a \textit{pseudolink} is an equivalence 
class of pseudolink diagrams under the equivalence relation generated by 
planar isotopy moves, the classical Reidemeister moves I, II and III, 
\[\includegraphics{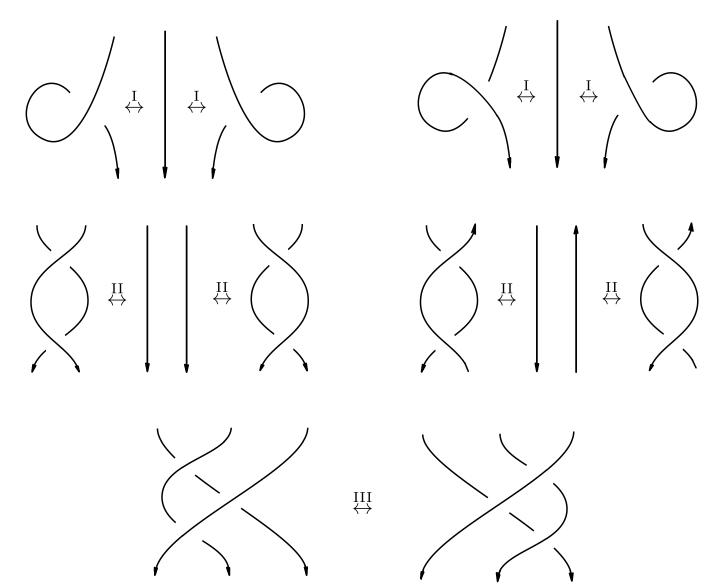}\]
and the moves PI, PI$'$, PII, PIII and PIII$'$ 
\[\includegraphics{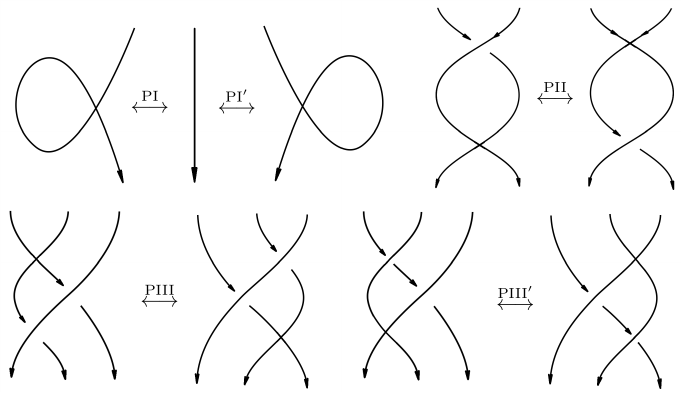}.\]
In particular, in \cite{NOS} it is shown that these moves form a generating 
set of oriented pseduoknot Reidemeister moves.
A \textit{singular link} is an equivalence class of singular link diagrams 
under the equivalence relation generated by planar isotopy, the usual 
Reidemeister moves I, II and III, and the moves PII, PIII and PIII$'$.
See \cite{BEHY,H,HHJJMR} for more about pseudoknot and singular knot
Reidemeister moves.

\begin{remark}
Singular knots and links may be regarded as the ``pseudo-framed case'' of 
pseudoknots and pseudolinks, where the ``pseudo-writhe'' or number of 
precrossings is preserved.
\end{remark}

As we will see, the algebraic conditions on our psybracket structure coming
from moves PI and PI$'$ are already implied by move PII so the 
invariants we define will be valid for both pseudolinks and singular links; 
on the other hand, they will not be able to distinguish singular links which 
differ only by PI and PI$'$ moves.

\section{\large\textbf{Psybrackets}}\label{PSY}

We begin with a definition.

\begin{definition}
Let $X$ be a set. A \textit{psybracket} structure on $X$ consists of 
two maps $\langle,,\rangle_c,\langle,,\rangle_p:X\times X\times X\to X$
such that
\begin{itemize}
\item[(i)] For all $a,b,c\in X$ there exist unique $x,y,z,u,v\in X$ 
such that
\[\begin{array}{rcll}
\langle a,b,x\rangle_c & = & c & (i.i)\\
\langle a,y,b\rangle_c & = & c & (i.ii)\\
\langle z,a,b\rangle_c & = & c & (i.iii)\\
\langle u,b,c\rangle_p & = & b & (i.iv) \\ 
\langle a,b,v\rangle_p & = & b & (i.v), 
\end{array}\]
\item[(ii)] 
For all $a,b,c\in X$ we have
\[
\langle a,\langle a,b,c\rangle_c,c\rangle_p =\langle a,\langle a,b,c\rangle_p,c\rangle_c 
\]
and
\item[(iii)] 
For all $a,b,c,d\in X$ we have
\[\begin{array}{rcll}
\langle \langle a,b,c\rangle_c,c,d\rangle_c  & = & 
\langle \langle  a,b,\langle b,c,d\rangle_c \rangle_c, \langle b,c,d\rangle_c, d\rangle_c & (iii.i)\\
& = & 
\langle \langle  a,b,\langle b,c,d\rangle_p \rangle_c, \langle b,c,d\rangle_p, d\rangle_c &(iii.ii)\\
\langle \langle a,b,c\rangle_p,c,d\rangle_c  & = & 
\langle \langle  a,b,\langle b,c,d\rangle_c \rangle_c, \langle b,c,d\rangle_c, d\rangle_p &(iii.iii)\\
\langle a,b,\langle b,c,d\rangle_c\rangle_c & = &
\langle a,\langle a,b,c\rangle_c,\langle \langle a,b,c\rangle_c,c,d\rangle_c\rangle_c &(iii.iv)\\
& = &
\langle a,\langle a,b,c\rangle_p,\langle \langle a,b,c\rangle_p,c,d\rangle_c\rangle_c &(iii.v)\\
\langle a,b,\langle b,c,d\rangle_p\rangle_c & = &
\langle a,\langle a,b,c\rangle_c,\langle \langle a,b,c\rangle_c,c,d\rangle_c\rangle_p &(iii.vi) 
\end{array}\]
\end{itemize}
\end{definition}

\begin{remark}
The set $X$ with the map $\langle,,\rangle_c$ forms a \textit{vertical 
tribracket} as described in \cite{NOO}, also known as a \textit{Niebrzydowski 
tribracket} or a \textit{knot-theoretic ternary quasigroup}. In particular, the
binary operations $\cdot_a, \cdot_a', \cdot_a'':X\times X\to X$ defined by
\[
x\cdot_a y= \langle a,x,y\rangle_c, \quad
x\cdot_a' y= \langle x,a,y\rangle_c,\quad \mathrm{and}\quad
x\cdot_a'' y= \langle x,y,a\rangle_c,
 \]
are quasigroup structures on $X$.
\end{remark}

The psybracket axioms are motivated by the following region coloring rules 
for pseudoknot diagrams:

\[\includegraphics{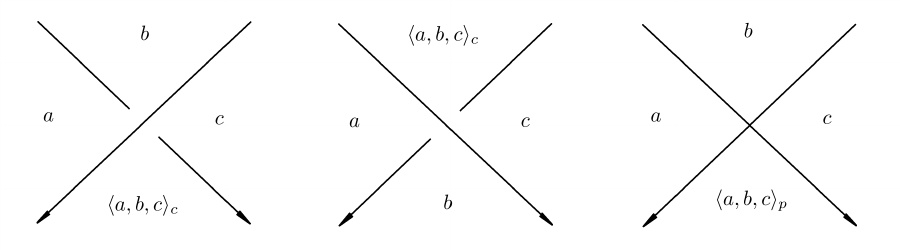}\]

An assignment of elements of a psybracket $X$ to each region in the complement 
of a pseudoknot diagram or singular knot diagram $D$ such that the above 
conditions are satisfied at every crossing is called a \textit{psybracket 
coloring} or an \textit{$X$-coloring} of $D$.
The psybracket axioms are the conditions needed to ensure that
for every psybracket coloring of a diagram before a move, there is a unique
coloring of the diagram after the move that agrees with the pre-move coloring 
outside the neighborhood of the move. 

Axioms (i.i)--(i.iv) are the conditions required by the Reidemeister II moves:
\[\scalebox{0.8}{\includegraphics{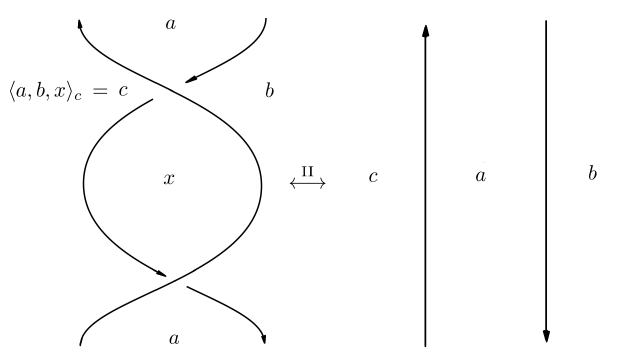}\quad \includegraphics{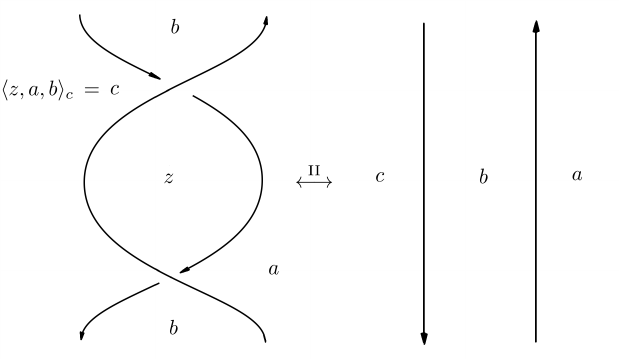}}\]
\[\scalebox{0.75}{\includegraphics{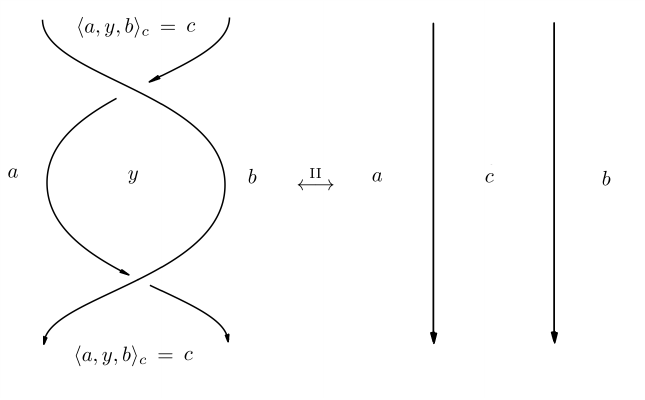} \quad \includegraphics{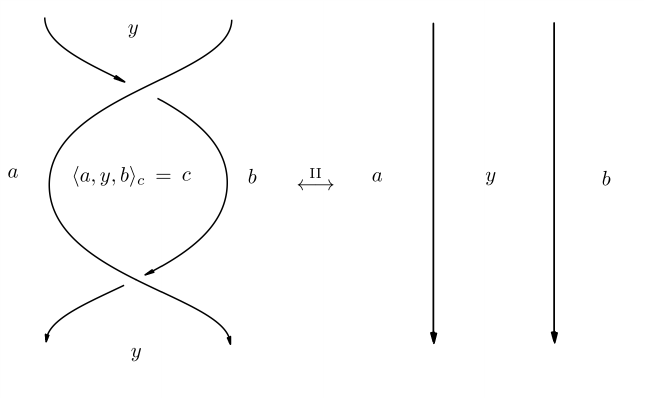}}\]

The conditions imposed by the classical Reidemeister I moves are special cases
of the Reidemeister II condtions. For example, in
\[\includegraphics{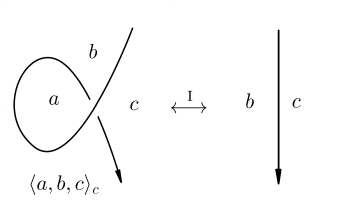}\]
we have the requirement that for every $b,c\in X$ there exists a unique $a\in X$
such that $\langle a,b,c\rangle _c=b$, but this is alread implied by 
left-invertibility (i.iii). The other classical Reidemeister I moves are 
similar.

Axioms (i.iv) and (i.v) are motivated by moves PI and PI$'$:
\[\includegraphics{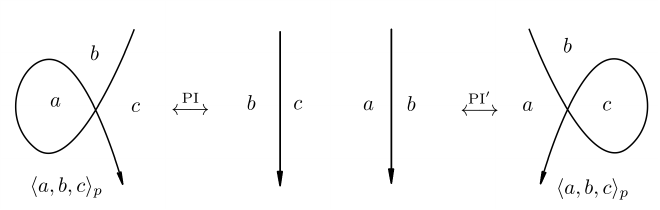}\]

Axiom (ii) is motivated by move PII:
\[\includegraphics{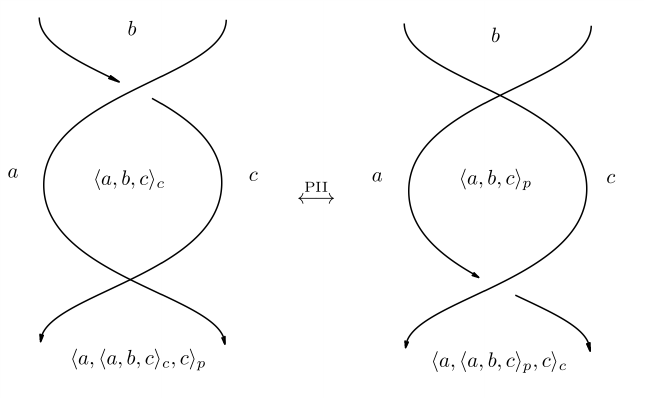}\]

Axioms (iii.i)--(iii.vi) are motivated by Reidemeister III,
PIII and PIII$'$:
\[\includegraphics{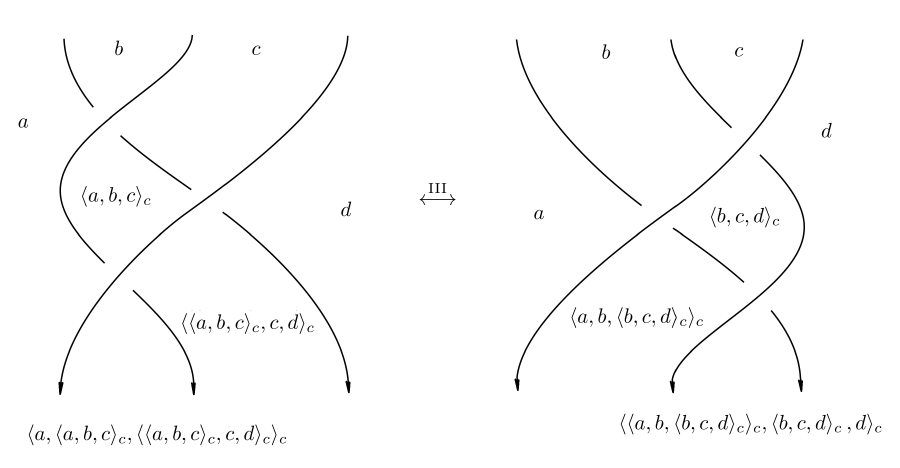}\]
\[\includegraphics{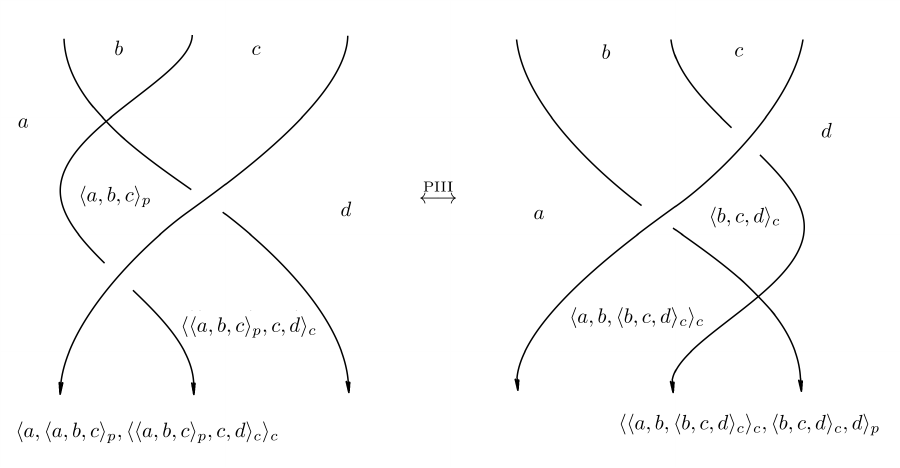}\]
\[\includegraphics{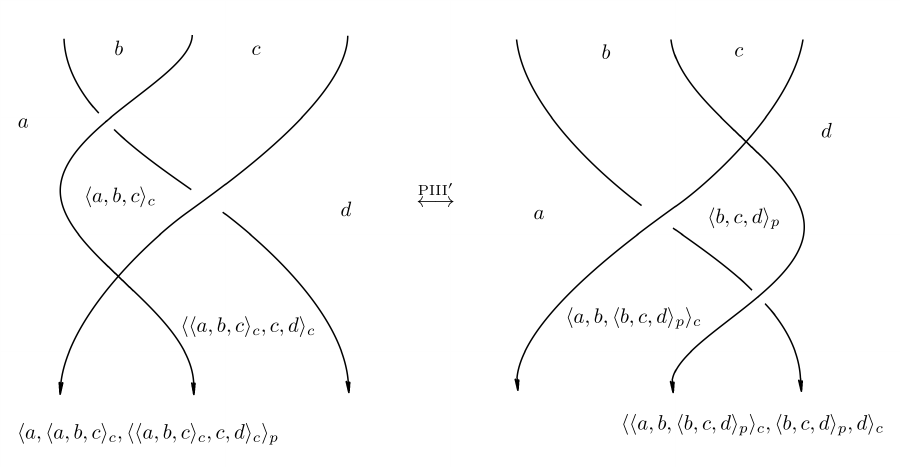}\]

Since this set of moves is a generating set of oriented pseudoknot Reidemeister
moves, we have:

\begin{proposition}\label{pm}
Let $D$ be an oriented pseduoknot diagram or singular knot diagram and $X$ a 
psybracket. Then for any $X$-coloring of $D$ before a psueduoknot or singular 
knot Reidemeister move, there is a unique $X$-coloring of the diagram after 
the move that agrees with the coloring before the move outside the 
neighborhood of the move.
\end{proposition}

\begin{definition}
Let $X$ and $Y$ be psybrackets. A map $f:X\to Y$ is a \textit{psybracket 
homomorphism} if for all $a,b,c\in X$ we have
\[\langle f(a),f(b),f(c)\rangle_c=f(\langle a,b,c\rangle_c) \ \mathrm{and}\
\langle f(a),f(b),f(c)\rangle_p=f(\langle a,b,c\rangle_p).
\]
A bijective psybracket homomorphism is an \textit{isomorphism}.
\end{definition}

\begin{example}\label{ex:pb1}
A Niebrzydowski tribracket can be given the structure of a psybracket 
by setting \[\langle a,b,c\rangle_p=\langle a,b,c\rangle_c.\] To see that
this definition satisfies the axioms, we need only note that replacing
the precrossings in the moves PI, PIII and PIII$'$ with positive crossings
results in
valid classical Reidemeister moves and in the same diagram on both sides 
of move PII. Similarly, setting
\[\langle a,b,c\rangle_p=d\]
where \[\langle a,d,c\rangle_c=b\]
yields a psybracket, as we can see by resolving the precrossings as
negative classical crossings.
\end{example}

\begin{example}
Let $G$ be a group. Then $G$ is a Niebrzydowski tribracket under the operation
\[\langle a,b,c\rangle =  ab^{-1}c\]
known as a \textit{Dehn tribracket}. Then the two psybracket structures
in Example \ref{ex:pb1} are 
\[\langle a,b,c\rangle_c = \langle a,b,c\rangle_p= ab^{-1}c\]
and
\[\langle a,b,c\rangle_c = ab^{-1}c,\quad  \langle a,b,c\rangle_p= cb^{-1}a\]
respectively; we call these the \textit{positive} and \textit{negative
Dehn psybracket} structures on $G$.
\end{example}

We can specify a psybracket structure on a finite set $X=\{1,2,\dots, n\}$
with a pair of \textit{operation 3-tensors},
i.e. lists of $n$ square matrices of size $n\times n$ such that the
entry in matrix $a$ row $b$ column $c$ is $\langle a,b,c\rangle_c$ or
$\langle a,b,c\rangle_p$.

\begin{example}
The Dehn tribracket structure on $\mathbb{Z}_4=\{1,2,3,4\}$ (where we use
$4$ for the class of zero) is
\[\left[
\left[\begin{array}{rrrr}
1 & 2 & 3 & 4 \\
4 & 1 & 2 & 3 \\
3 & 4 & 1 & 2 \\
2 & 3 & 4 & 1
\end{array}\right],
\left[\begin{array}{rrrr}
2 & 3 & 4 & 1 \\
1 & 2 & 3 & 4 \\
4 & 1 & 2 & 3 \\
3 & 4 & 1 & 2 
\end{array}\right],
\left[\begin{array}{rrrr}
3 & 4 & 1 & 2 \\
2 & 3 & 4 & 1 \\
1 & 2 & 3 & 4 \\
4 & 1 & 2 & 3 
\end{array}\right],
\left[\begin{array}{rrrr}
4 & 1 & 2 & 3 \\
3 & 4 & 1 & 2 \\
2 & 3 & 4 & 1 \\
1 & 2 & 3 & 4
\end{array}\right]
\right].\]
The positive and negative Dehn psybracket structures on $\mathbb{Z}_4$
both have this same operation 3-tensor for 
$\langle,,\rangle_p$ and for $\langle,,\rangle_c$.
\end{example}

\begin{example}
Let $X=\{1,2,3\}$. Using \textit{python} code, we compute that there are six
isomorphism classes of psybracket structures on $X$. Representatives
of each class are specified by the operation 3-tensors below:

\[\left[\left[\left[\begin{array}{rrr}
1 & 2& 3\\ 2& 3& 1\\ 3& 1& 2
\end{array}\right],\left[\begin{array}{rrr}
3 & 1& 2\\ 1& 2& 3\\ 2& 3& 1
\end{array}\right],\left[\begin{array}{rrr}
2 & 3& 1\\ 3& 1& 2\\ 1& 2& 3
\end{array}\right]\right]_c, 
\left[\left[\begin{array}{rrr}
1 & 2& 3\\ 2& 3& 1\\ 3& 1& 2
\end{array}\right],\left[\begin{array}{rrr}
3 & 1& 2\\ 1& 2& 3\\ 2& 3& 1
\end{array}\right],\left[\begin{array}{rrr}
2 & 3& 1\\ 3& 1& 2\\ 1& 2& 3
\end{array}\right]\right]_p\right]
\]
\[\left[\left[\left[\begin{array}{rrr}
1 & 2& 3\\ 2& 3& 1\\ 3& 1& 2
\end{array}\right],\left[\begin{array}{rrr}
3 & 1& 2\\ 1& 2& 3\\ 2& 3& 1
\end{array}\right],\left[\begin{array}{rrr}
2 & 3& 1\\ 3& 1& 2\\ 1& 2& 3
\end{array}\right]\right]_c, 
\left[\left[\begin{array}{rrr}
1 & 3 & 2 \\2 & 1 & 3 \\3 & 2 & 1
\end{array}\right],\left[\begin{array}{rrr}
2 & 1 & 3 \\3 & 2 & 1 \\1 & 3 & 2
\end{array}\right],\left[\begin{array}{rrr}
3 & 2 & 1 \\1 & 3 & 2 \\2 & 1 & 3
\end{array}\right]\right]_p\right]
\]
\[\left[\left[\left[\begin{array}{rrr}
2& 3& 1 \\3& 1& 2 \\1& 2& 3
\end{array}\right],\left[\begin{array}{rrr}
1& 2& 3\\2& 3& 1 \\3& 1& 2
\end{array}\right],\left[\begin{array}{rrr}
3& 1& 2 \\1& 2& 3 \\2& 3& 1
\end{array}\right]\right]_c,
\left[\left[\begin{array}{rrr}
2& 3& 1 \\3& 1& 2 \\1& 2& 3
\end{array}\right],\left[\begin{array}{rrr}
1& 2& 3 \\2& 3& 1 \\3& 1& 2
\end{array}\right],\left[\begin{array}{rrr}
3& 1& 2 \\1& 2& 3 \\2& 3& 1
\end{array}\right]\right]_p\right]
\]
\[\left[\left[\left[\begin{array}{rrr}
2 & 3 & 1\\ 3 & 1 & 2\\ 1 & 2 & 3
\end{array}\right],\left[\begin{array}{rrr}
1 & 2 & 3\\ 2 & 3 & 1\\ 3 & 1 & 2
\end{array}\right],\left[\begin{array}{rrr}
3 & 1 & 2\\ 1 & 2 & 3\\ 2 & 3 & 1
\end{array}\right]\right]_c,
\left[\left[\begin{array}{rrr}
3 & 2 & 1\\ 1 & 3 & 2\\ 2 & 1 & 3
\end{array}\right],\left[\begin{array}{rrr}
1 & 3 & 2\\ 2 & 1 & 3\\ 3 & 2 & 1
\end{array}\right],\left[\begin{array}{rrr}
2 & 1 & 3\\ 3 & 2 & 1\\ 1 & 3 & 2
\end{array}\right]\right]_p\right]
\]
\[\left[\left[\left[\begin{array}{rrr}
1 & 3 & 2 \\3 & 2 & 1 \\2 & 1 & 3
\end{array}\right],\left[\begin{array}{rrr}
3 & 2 & 1 \\2 & 1 & 3 \\1 & 3 & 2
\end{array}\right],\left[\begin{array}{rrr}
2 & 1 & 3 \\1 & 3 & 2 \\3 & 2 & 1
\end{array}\right]\right]_c,
\left[\left[\begin{array}{rrr}
1 & 2 & 3 \\1 & 2 & 3 \\1 & 2 & 3
\end{array}\right],\left[\begin{array}{rrr}
2 & 3 & 1 \\2 & 3 & 1 \\2 & 3 & 1
\end{array}\right],\left[\begin{array}{rrr}
3 & 1 & 2 \\3 & 1 & 2 \\3 & 1 & 2
\end{array}\right]\right]_p\right]
\]
\[\left[\left[\left[\begin{array}{rrr}
1 & 3 & 2 \\ 3 & 2 & 1 \\ 2 & 1 & 3
\end{array}\right],\left[\begin{array}{rrr}
3 & 2 & 1 \\ 2 & 1 & 3 \\ 1 & 3 & 2
\end{array}\right],\left[\begin{array}{rrr}
2 & 1 & 3 \\ 1 & 3 & 2 \\ 3 & 2 & 1
\end{array}\right]\right]_c,
\left[\left[\begin{array}{rrr}
1 & 3 & 2 \\ 3 & 2 & 1 \\ 2 & 1 & 3
\end{array}\right],\left[\begin{array}{rrr}
3 & 2 & 1 \\ 2 & 1 & 3 \\ 1 & 3 & 2
\end{array}\right],\left[\begin{array}{rrr}
2 & 1 & 3 \\ 1 & 3 & 2 \\ 3 & 2 & 1
\end{array}\right]\right]_p\right]
\]
\end{example}

\section{\large\textbf{Psybracket Counting Invariants}}\label{PCI}

The motivation for the psybracket definition is to define counting 
invariants of pseudoknots and singular knots. 
As with previous knot coloring structures, by Proposition \ref{pm} we have the 
following result:

\begin{theorem}
Let $X$ be a finite psybracket, $K$ an oriented pseudoknot or singular knot
diagram, and $\mathcal{C}(K,X)$ the set of $X$-colorings of $K$. Then the 
number of $X$-colorings of $K$, 
\[\Phi_X^{\mathbb{Z}}(K)=|\mathcal{C}(K,X)|\]
is an integer-valued invariant of pseudoknots and singular knots we will 
call the \textit{psybracket counting invariant}. 
\end{theorem}

\begin{example}
Let $X$ be the psybracket specified by 
\[\left[\left[\left[\begin{array}{rrr}
1 & 3 & 2 \\ 3 & 2 & 1 \\ 2 & 1 & 3
\end{array}\right],\left[\begin{array}{rrr}
3 & 2 & 1 \\ 2 & 1 & 3 \\ 1 & 3 & 2
\end{array}\right],\left[\begin{array}{rrr}
2 & 1 & 3 \\ 1 & 3 & 2 \\ 3 & 2 & 1
\end{array}\right]\right]_c,
\left[\left[\begin{array}{rrr}
1 & 3 & 2 \\ 3 & 2 & 1 \\ 2 & 1 & 3
\end{array}\right],\left[\begin{array}{rrr}
3 & 2 & 1 \\ 2 & 1 & 3 \\ 1 & 3 & 2
\end{array}\right],\left[\begin{array}{rrr}
2 & 1 & 3 \\ 1 & 3 & 2 \\ 3 & 2 & 1
\end{array}\right]\right]_p\right]
\]
and consider the pseudolink 
\[\includegraphics{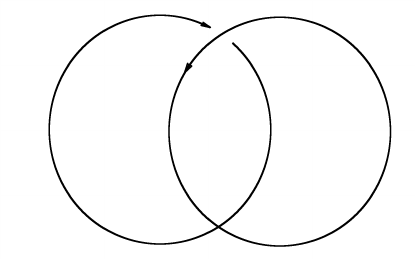}\]
Then for example the assignment of elements in $\{1,2,3\}$ to the
regions given by
\[\includegraphics{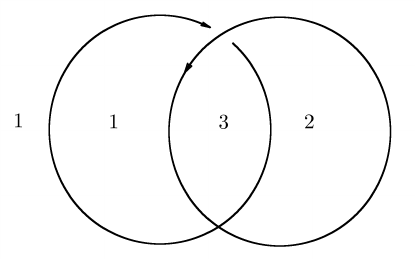}\]
is an $X$-coloring, since we have
\[\langle 1,1,2\rangle_c=3 \quad \mathrm{and}\quad
\langle 1,3,2\rangle_p=1.  \]
\end{example}

\begin{example}
Consider the pseudoknot $3_1.3$ below. 
\[\includegraphics{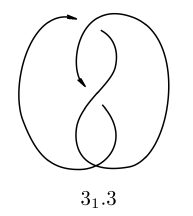}\]
It has two positive crossings and one
precrossing; resolving the precrossing one way yields a trefoil and the other
way yields an unknot. We can distinguish this pseudoknot from both the
trefoil and the unknot using psybracket counting invariants. Specifically,
the psybracket 
\[X_1=\left[\left[\left[\begin{array}{rrr}
1 & 2& 3\\ 2& 3& 1\\ 3& 1& 2
\end{array}\right],\left[\begin{array}{rrr}
3 & 1& 2\\ 1& 2& 3\\ 2& 3& 1
\end{array}\right],\left[\begin{array}{rrr}
2 & 3& 1\\ 3& 1& 2\\ 1& 2& 3
\end{array}\right]\right]_c, 
\left[\left[\begin{array}{rrr}
1 & 2& 3\\ 2& 3& 1\\ 3& 1& 2
\end{array}\right],\left[\begin{array}{rrr}
3 & 1& 2\\ 1& 2& 3\\ 2& 3& 1
\end{array}\right],\left[\begin{array}{rrr}
2 & 3& 1\\ 3& 1& 2\\ 1& 2& 3
\end{array}\right]\right]_p\right]
\]
gives us counting invariant values 
$\Phi_{X_1}^{\mathbb{Z}}(3_1.3)=27$ and $\Phi_{X_1}^{\mathbb{Z}}(\mathrm{Unknot})=9$ 
while the psybracket
\[X_2=\left[\left[\left[\begin{array}{rrr}
1 & 2& 3\\ 2& 3& 1\\ 3& 1& 2
\end{array}\right],\left[\begin{array}{rrr}
3 & 1& 2\\ 1& 2& 3\\ 2& 3& 1
\end{array}\right],\left[\begin{array}{rrr}
2 & 3& 1\\ 3& 1& 2\\ 1& 2& 3
\end{array}\right]\right]_c, 
\left[\left[\begin{array}{rrr}
1 & 3 & 2 \\2 & 1 & 3 \\3 & 2 & 1
\end{array}\right],\left[\begin{array}{rrr}
2 & 1 & 3 \\3 & 2 & 1 \\1 & 3 & 2
\end{array}\right],\left[\begin{array}{rrr}
3 & 2 & 1 \\1 & 3 & 2 \\2 & 1 & 3
\end{array}\right]\right]_p\right]
\] 
gives us counting invariant values 
$\Phi_{X_2}^{\mathbb{Z}}(3_1.3)=9$ and $\Phi_{X_2}^{\mathbb{Z}}(3_1)=27$. 
\end{example}

\begin{example}
	We computed the coloring invariant $\Phi_{X}^{\mathbb{Z}}$ of a choice of orientation for the 2-bouquet graphs in \cite{O} using the psybracket $X$ with the operation matrix
	
	\[\left[\left[\left[\begin{array}{rrr}
	1 & 2& 3\\ 2& 3& 1\\ 3& 1& 2
	\end{array}\right],\left[\begin{array}{rrr}
	3 & 1& 2\\ 1& 2& 3\\ 2& 3& 1
	\end{array}\right],\left[\begin{array}{rrr}
	2 & 3& 1\\ 3& 1& 2\\ 1& 2& 3
	\end{array}\right]\right]_c, 
	\left[\left[\begin{array}{rrr}
	1 & 2& 3\\ 2& 3& 1\\ 3& 1& 2
	\end{array}\right],\left[\begin{array}{rrr}
	3 & 1& 2\\ 1& 2& 3\\ 2& 3& 1
	\end{array}\right],\left[\begin{array}{rrr}
	2 & 3& 1\\ 3& 1& 2\\ 1& 2& 3
	\end{array}\right]\right]_p\right].
	\]
	
	The results are collected in the table.
	\begin{center}
		
		\begin{tabular}{ c|c }
			$L$& $\Phi_{X}^{\mathbb{Z}}(L)$ \\
			\hline
$0_1^k, 5^k_8,6^k_4, 6^k_5,1^l_1,6^l_{10}$ &$ 9$\\ \hline
$6^k_6, 6^k_7,5^l_1,6^l_7,6^l_{12}$  &$ 27 $\\ \hline
$2^k_1,3^k_1,5^k_6, 5^k_7, 6^k_{16},6^k_{17},6^k_{19},4^l_1,5^l_3,6^l_1, 6^l_3$ &$ 81 $\\ \hline
$5^k_2,5^k_4,5^k_5,6^k_{15},6^k_{19},3^l_1,6^l_2,6^l_4, 6^l_5,6^l_8$ &$ 243 $\\ \hline
$4^k_1, 4^k_2,4^k_3, 5^k_1, 5^k_3,6^k_{11},6^k_{13}, 6^k_{14},6^k_{18},6^l_6,6^l_9$ & $729 $\\ \hline
$6^k_{12},5^l_2, 6^l_{11}$ & $2187 $\\ \hline
$6^k_1,6^k_2, 6^k_3,6^k_8, 6^k_9,6^k_{10}$  & $6561 $ 
		\end{tabular}
	\end{center}
\end{example}

\begin{example}
	We computed the coloring invariant $\Phi_{X}^{\mathbb{Z}}$ of a choice of orientation for the pseudoknots up to 5 crossings in the pseudoknot tables in \cite{HHJJMR} using the psybracket $X$ with the operation matrix
	
	\[\left[\left[\left[\begin{array}{rrr}
	1 & 2& 3\\ 2& 3& 1\\ 3& 1& 2
	\end{array}\right],\left[\begin{array}{rrr}
	3 & 1& 2\\ 1& 2& 3\\ 2& 3& 1
	\end{array}\right],\left[\begin{array}{rrr}
	2 & 3& 1\\ 3& 1& 2\\ 1& 2& 3
	\end{array}\right]\right]_c, 
	\left[\left[\begin{array}{rrr}
	1 & 2& 3\\ 2& 3& 1\\ 3& 1& 2
	\end{array}\right],\left[\begin{array}{rrr}
	3 & 1& 2\\ 1& 2& 3\\ 2& 3& 1
	\end{array}\right],\left[\begin{array}{rrr}
	2 & 3& 1\\ 3& 1& 2\\ 1& 2& 3
	\end{array}\right]\right]_p\right].
	\]

	The results are collected in the table.
	\begin{center}
		
		\begin{tabular}{ c|c }
			$L$& $\Phi_{X}^{\mathbb{Z}}(L)$ \\
			\hline
			$5_1.1, 5_1.2, 5_1.3,5_1.4,5_1.5,5_2.1$ & 9 \\ 
			\hline
			$3_1.1,3_1.2,3_1.3,4_1.1,4_1.2,4_1.3,4_1.4,4_1.5,5_2.2$ &27 \\ 
			\hline
			$5_2.3,5_2.4,5_2.6$ & 81\\
			\hline
			$5_2.5,5_2.7,5_2.9$ & 243\\
			\hline
			$5_2.8,5_2.10$ & 729\\
		\end{tabular}
	\end{center}
\end{example}

\begin{example}
	
We computed the coloring invariant $\Phi_{X}^{\mathbb{Z}}$ of a choice of orientation for the pseudoknots up to 5 crossings in the pseudoknot tables using the psybracket $X$ with the operation matrix

\[\left[\left[\left[\begin{array}{rrr}
1 & 2& 3\\ 2& 3& 1\\ 3& 1& 2
\end{array}\right],\left[\begin{array}{rrr}
3 & 1& 2\\ 1& 2& 3\\ 2& 3& 1
\end{array}\right],\left[\begin{array}{rrr}
2 & 3& 1\\ 3& 1& 2\\ 1& 2& 3
\end{array}\right]\right]_c, 
\left[\left[\begin{array}{rrr}
1 & 3 & 2 \\2 & 1 & 3 \\3 & 2 & 1
\end{array}\right],\left[\begin{array}{rrr}
2 & 1 & 3 \\3 & 2 & 1 \\1 & 3 & 2
\end{array}\right],\left[\begin{array}{rrr}
3 & 2 & 1 \\1 & 3 & 2 \\2 & 1 & 3
\end{array}\right]\right]_p\right].
\]

The results are collected in the table.
\begin{center}
	
	\begin{tabular}{ c|c }
		$L$& $\Phi_{X}^{\mathbb{Z}}(L)$ \\
		\hline
		$3_1.2,3_1.3,4_1.2,4_1.3,5_1.1,5_1.3,5_1.4,5_2.1$ & 9 \\ 
		\hline
		$3_1.1,4_1.1,4_1.4,4_1.5,5_1.2,5_1.5,5_2.2$ &27 \\ 
		\hline
		$5_2.3,5_2.4,5_2.6$ & 81\\
		\hline
		$5_2.5,5_2.7,5_2.9$ & 243\\
		\hline
		$5_2.8,5_2.10$ & 729\\
	\end{tabular}
\end{center}

\end{example}

\section{\large\textbf{Questions}}\label{Q}

In this paper we have only initiated the study of the new topic of 
psybrackets and their pseudo/singular knot and link invariants. There
are many interesting questions to be explored in this area; we suggest
a few of them here.

\begin{itemize}
\item As with counting invariants arising from other structures, many types
of \textit{enhancements} are possible. Applying a historically successful 
strategy, we ask what invariants of psybracket-colored pseudoknots are possible.
Ideas might include cocycle enhancements analogous to those in \cite{NOO},
skein enhancements like those in \cite{ANR}, module enhancements like those in
\cite{NNS} and many more.
\item What is the structure of psybrackets? What kinds of products, 
decompositions, functors to and from other algebraic categories are possible?
\item What generalizations are possible to the cases of pseudo/singular 
trivalent graphs and handlebody knots or to the virtual and twisted virtual cases?
\end{itemize}

\bibliography{jk-sj-sn}{}

\begin{thebibliography}{10}

\bibitem{ANR}
L.~Aggarwal, S.~Nelson, and P.~Rivera.
\newblock Quantum enhancements via tribracket brackets.
\newblock {\em arXiv:1907.03011}, 2019.

\bibitem{BEHY}
K.~Bataineh, M.~Elhamdadi, M.~Hajij, and W.~Youmans.
\newblock Generating sets of reidemeister moves of oriented singular links and
  quandles.
\newblock {\em arXiv:1702.01150}, 2017.

\bibitem{CNS}
W.~Choi, D.~Needell, and S.~Nelson.
\newblock Boltzmann enhancements of biquasile counting invariants.
\newblock {\em J. Knot Theory Ramifications}, 27(14):1850068, 12, 2018.

\bibitem{CDDK}
P.~Clote, S.~Dobrev, I.~Dotu, E.~Kranakis, D.~Krizanc, and J.~Urrutia.
\newblock On the page number of {RNA} secondary structures with pseudoknots.
\newblock {\em J. Math. Biol.}, 65(6-7):1337--1357, 2012.

\bibitem{EP}
P.~A. Evans.
\newblock Finding common {RNA} pseudoknot structures in polynomial time.
\newblock {\em J. Discrete Algorithms}, 9(4):335--343, 2011.

\bibitem{GNT}
P.~Graves, S.~Nelson, and S.~Tamagawa.
\newblock Niebrzydowski algebras and trivalent spatial graphs.
\newblock {\em Internat. J. Math.}, 29(14):1850102, 16, 2018.

\bibitem{H}
R.~Hanaki.
\newblock Pseudo diagrams of knots, links and spatial graphs.
\newblock {\em Osaka J. Math.}, 47(3):863--883, 2010.

\bibitem{HHJJMR}
A.~Henrich, R.~Hoberg, S.~Jablan, L.~Johnson, E.~Minten, and L.~Radovi\'c.
\newblock The theory of pseudoknots.
\newblock {\em J. Knot Theory Ramifications}, 22(7):1350032, 21, 2013.

\bibitem{HJ}
A.~Henrich and S.~Jablan.
\newblock On the coloring of pseudoknots.
\newblock {\em J. Knot Theory Ramifications}, 23(12):1450061, 22, 2014.

\bibitem{KN}
J.~Kim and S.~Nelson.
\newblock Biquasile colorings of oriented surface-links.
\newblock {\em Topology Appl.}, 236:64--76, 2018.

\bibitem{LR}
T.~J.~X. Li and C.~M. Reidys.
\newblock Statistics of topological {RNA} structures.
\newblock {\em J. Math. Biol.}, 74(7):1793--1821, 2017.

\bibitem{DSN2}
D.~Needell and S.~Nelson.
\newblock Biquasiles and dual graph diagrams.
\newblock {\em J. Knot Theory Ramifications}, 26(8):1750048, 18, 2017.

\bibitem{NNS}
D.~Needell, S.~Nelson, and Y.~Shi.
\newblock Tribracket {M}odules.
\newblock {\em Internat. J. Math.}, 31(4):2050028, 13, 2020.

\bibitem{NOR}
S.~Nelson, M.~E. Orrison, and V.~Rivera.
\newblock Quantum enhancements and biquandle brackets.
\newblock {\em J. Knot Theory Ramifications}, 26(5):1750034, 24, 2017.

\bibitem{NOO}
S.~Nelson, K.~Oshiro, and N.~Oyamaguchi.
\newblock Local biquandles and {N}iebrzydowski's tribracket theory.
\newblock {\em Topology Appl.}, 258:474--512, 2019.

\bibitem{NOS}
S.~Nelson, N.~Oyamaguchi, and R.~Sazdanovic.
\newblock Psyquandles, singular knots and pseduoknots.
\newblock {\em Tokyo J. Math.}, 42(2):405--429, 2019.

\bibitem{NP}
S.~Nelson and E.~Pauletich.
\newblock Multi-tribrackets.
\newblock {\em J. Knot Theory Ramifications}, 28(12):1950075, 16, 2019.

\bibitem{NP2}
S.~Nelson and S.~Pico.
\newblock Virtual tribrackets.
\newblock {\em J. Knot Theory Ramifications}, 28(4):1950026, 12, 2019.

\bibitem{MN}
M.~Niebrzydowski.
\newblock On some ternary operations in knot theory.
\newblock {\em Fund. Math.}, 225(1):259--276, 2014.

\bibitem{MN2}
M.~Niebrzydowski, A.~Pilitowska, and A.~Zamojska-Dzienio.
\newblock Knot-theoretic ternary groups.
\newblock {\em Fund. Math.}, 247(3):299--320, 2019.

\bibitem{O}
N.~Oyamaguchi.
\newblock Enumeration of spatial 2-bouquet graphs up to flat vertex isotopy.
\newblock {\em Topology Appl.}, 196(part B):805--814, 2015.

\bibitem{QR}
J.~Qin and C.~M. Reidys.
\newblock On topological {RNA} interaction structures.
\newblock {\em J. Comput. Biol.}, 20(7):495--513, 2013.

\end{thebibliography}
\bibliographystyle{abbrv}

\bigskip

\noindent
\textsc{Department of Mathematics, \\
Pusan National University, \\
Busan 46241, Republic of Korea
}

\medskip

\noindent
\textsc{Department of Mathematical Sciences \\
Claremont McKenna College \\
850 Columbia Ave. \\
Claremont, CA 91711}

\end{document}